\documentclass[reqno]{amsart}
\usepackage{amssymb,amsmath,amsthm,amsfonts,fancyhdr}
\usepackage{indentfirst}
\usepackage{mathrsfs}
\usepackage{setspace}
\usepackage{color}
\usepackage{cite}
\usepackage{bm}
\numberwithin{equation}{section}

\newtheorem{thm}{Theorem}[section]
\newtheorem{defn}[thm]{Defnition}
\newtheorem{lemma}[thm]{Lemma}
\newtheorem{cor}[thm]{Corollary}
\newtheorem{re}[thm]{Remark}

  \newcommand{\dist}{\mbox{dist}}
  \newcommand{\diam}{{\mbox{diam}}}

  \numberwithin{equation}{section}
  \numberwithin{figure}{section}

\begin{document}

\title[${{W}^{2,\delta}}$ Estimates for Fully Nonlinear Elliptic inequalities]{$\bm{{W}^{2,\delta}}$ Estimates for Solution Sets of Fully Nonlinear Elliptic inequalities on $\bm{C^{1,\alpha}}$ Domains}

\author{Dongsheng Li}
\author{Xuemei Li}

\address{School of Mathematics and Statistics, Xi'an Jiaotong University, Xi'an, P.R.China 710049.}
\address{School of Mathematics and Statistics, Xi'an Jiaotong University, Xi'an, P.R.China 710049.}

\email{lidsh@mail.xjtu.edu.cn}
\email{pwflxm@stu.xjtu.edu.cn}

\begin{abstract}
In this paper,  we establish boundary $W^{2,\delta}$ estimates for $u\in S(\lambda,\Lambda,f)$  on $C^{1,\alpha}$ domains with $f\in L^p$ as $n<p<\infty$ and $C^{1,\alpha}$ boundary values.
Instead of straightening out the boundary, our main idea is to obtain boundary $W^{2,\delta}$ estimates from interior $W^{2,\delta_0}$ estimates and Whitney decomposition for some $\delta\leq \delta_0$.
\end{abstract}



\maketitle

\section{Introduction}
\label{}
In this paper, we establish boundary $W^{2,\delta}$ estimates for $u\in S(\lambda,\Lambda,f)$ on $C^{1,\alpha}$ domains. Some notions and notations concerning  $S(\lambda,\Lambda,f)$ are listed as follows. We refer to \cite{CC} and \cite{CCKS} for more details.

\begin{defn}
Let $0<\lambda\leq\Lambda<\infty$. For $M\in \mathcal{S}(n)$ (the space of real $n\times n$ symmetric matrices), the Pucci's extremal operators $\mathcal{M}^{-}$ and $\mathcal{M}^{+}$ are defined as follows:
$$\mathcal{M}^{-}(M,\lambda,\Lambda)=\lambda\sum\limits_{e_i>0}e_i+\Lambda\sum\limits_{e_i<0}e_i,$$
$$\mathcal{M}^{+}(M,\lambda,\Lambda)=\Lambda\sum\limits_{e_i>0}e_i+\lambda\sum\limits_{e_i<0}e_i,$$
where $e_i=e_i(M)$ are the eigenvalues of $M$.
\end{defn}

\begin{defn}
Let $p>n/2$ and $f\in L^p_{loc}(\Omega)$. A function $u\in C(\Omega)$ satisfies  $\mathcal{M}^{-}(D^2u,\lambda,\Lambda)\leq f(x)$ (resp. $\mathcal{M}^{+}(D^2u,\lambda,\Lambda)\geq f(x)$) in $\Omega$ in the viscosity sense
if $x_0\in \Omega$, $\varphi\in W^{2,p}_{loc}(\Omega)$ and $u-\varphi$ has a local minimum (resp. maximum) at $x_0$ imply
$$ess\liminf\limits_{x\rightarrow x_0}\mathcal\{{M}^{-}(D^2 \varphi,\lambda,\Lambda)-f(x)\}\leq 0,$$
$$(resp.\ ess\limsup\limits_{x\rightarrow x_0}\mathcal\{{M}^{+}(D^2 \varphi,\lambda,\Lambda)-f(x)\}\geq 0).$$
\end{defn}

\begin{defn}
Let $p>n/2$ and $f\in L^p_{loc}(\Omega)$. We denote by $S(\lambda,\Lambda,f)$ the space of continuous functions $u$ in $\Omega$ such that
$$\mathcal{M}^{-}(D^2u,\lambda,\Lambda)\leq f(x)\leq \mathcal{M}^{+}(D^2u,\lambda,\Lambda)$$
in the viscosity sense in $\Omega$.
\end{defn}

Our main result is as follows.


\bigskip

\begin{thm}\label{t1.1}
 Let $0<\alpha< 1$ and $\Omega\subset\mathbb{R}^n$ be a bounded domain with $C^{1,\alpha}$ boundary portion $T\subset\partial\Omega$. Suppose that  $u$ satisfies
\begin{equation}
\left\{
\begin{aligned}\label{1.1} u&\in S(\lambda,\Lambda,f)   &\mathrm{in} \ \   &\Omega;\\
u&=g   &\mathrm{on}  \ \  &T,
\end{aligned}
\right.
\end{equation}
where
$$f\in  L^{p}(\Omega) \ \ \mathrm{for}\ \ n<p<\infty\ \ \mathrm{and}\ \  g\in C^{1,\alpha}(T).$$
Then there exists a positive constant $\delta$ depending on $n$, $\lambda$, $\Lambda$, $\alpha$ and $p$ such that for any domain $\Omega'\subset\subset\Omega\cup T$,
we have
\begin{equation}\label{1.2}
||u ||_{ W^{2,\delta}(\Omega')}\leq C\left( ||u ||_{L^{\infty}(\Omega)}+ ||f||_{L^{p}(\Omega)}+||g||_{C^{1,\alpha}(T)}\right),
\end{equation}
where $C$ depends on $n$, $\lambda$, $\Lambda$, $\alpha$, $\delta$, $p$, $T$, $\Omega'$ and $\Omega$.
\end{thm}

\bigskip

\begin{re}\label{r1.1}
$(i)$ When $T=\partial\Omega$ in Theorem \ref{t1.1}, we may obtain a global $W^{2,\delta}(\Omega)$ estimate on the bounded $C^{1,\alpha}$ domain.

(ii) Interior $W^{2,\delta_0}$ estimates and boundary H${\ddot{o}}$lder gradient estimates are used in the proof of Theorem \ref{t1.1}, where $\delta_0$ depends on $n,\lambda$ and $\Lambda$.
Interior $W^{2,\delta_0}$ estimates need $p>n-\epsilon$ with $\epsilon$ depending on $n,\lambda$ and $\Lambda$. (See \cite{E} for more details.) Boundary H${\ddot{o}}$lder gradient estimates need $p>n$. (See the following Corollary \ref{l3.3} for more details.)
\end{re}

\bigskip

$W^{2,p}$ regularity for fully nonlinear elliptic equations has been studied extensively during recent decades. For $u\in S(\lambda,\Lambda,f)$ and $f\in L^n$, Caffarelli \cite{CC} proved interior $W^{2,\delta_0}$ estimates with $\delta_0>0$ depending on $n,\lambda$ and $\Lambda$, which was first discovered by Lin \cite{L} for linear elliptic equations in nondivergence form. Then, interior $W^{2,p}$ estimates for fully nonlinear elliptic equations
$F(D^2 u,x)=f(x)$ with $f\in L^p$ and $n<p<\infty$ were established under the following two
 assumptions:  The homogeneous equations with constants coefficients $F(D^2 w,x_0)$ have  interior  $C^{1,1}$ estimates for any $x_0\in B_1$ and the oscillation of $F(M,x)$ in $x$ is small in $L^{n}$ sense. (See \cite{CC}.) Escauriaza\cite{E} extended the above estimates to the optimal range of exponents $n-\epsilon<p<\infty$ with $\epsilon$ depending on $n,\lambda$ and $\Lambda$. In \cite{DZ}, Li and Zhang relaxed the Caffarelli's assumptions of $F$ to: $F(D^2 w,x_0)$ have interior $W^{2,p_1}$ estimates for any $x_0\in B_1$ and the oscillation of $F(M,x)$ in $x$ is small in $L^{p_2}$ sense, where $n-\epsilon<p_0<p<p_1\leq\infty$ and $1/p_1+1/p_2=1/p_0$.  Winter \cite{W} extended Caffarelli's interior $W^{2,p}$ regularity to boundary
and established boundary $W^{2,p}$ estimates on flat domains and then on $C^{1,1}$ domains
alongside a flattening argument.
For oblique boundary value problem for $F(D^2 u,x)=f(x)$, based on pointwise boundary $C^{2,\alpha}$ regularity in \cite{LZ},
Byun and Han \cite{B} derived
the boundary $W^{2,p}$ estimates on flat domains and then on $C^{2,\alpha}$ domains by straightening out boundary.

The technique of straightening out boundary is widely applied to deal with boundary regularity. However, it is no longer applicable for the $C^{1,\alpha}$ domains. The reason is as follows: If the technique of straightening out boundary is used, then
for any $x_0\in \partial\Omega$, one should find a neighborhood $\mathcal N$  of $x_0$ and a diffeomorphism $\psi$ from $\mathcal N$ onto $B_1$ such that
$\psi(\mathcal N\cap \Omega)=B_1^+$ and $\psi(\mathcal N\cap \partial\Omega)= B_1\cap\{x^n=0\}$.
Writing $y=\psi(x),\tilde u(y)=u(x)$ and $\tilde f(y)=f(x)$ for $x\in \mathcal N\cap \Omega$ and $y\in B_1^+$,
we calculate
$$D^2u=D\psi^T (D^2\tilde u\circ\psi) D\psi
+((D\tilde u\circ\psi) D_{ij}\psi)_{1\leq i,j\leq n}.$$
Then
$\psi$ need to be $C^{1,1}$ in view of the second term.

In this paper, instead of straightening out the boundary, we establish boundary $W^{2,\delta}$ estimates for $u\in S(\lambda,\Lambda,f)$ on $C^{1,\alpha}$ domains from interior $W^{2,\delta_0}$ estimates and Whitney decomposition.
 Whitney decomposition is an effective tool for obtaining boundary estimates from interior estimates. For instance, Cao, Li and Wang \cite{CL} utilized it to prove the optimal weighted $W^{2,p}$  estimates for elliptic equations with non-compatible conditions; Li, Li and Zhang \cite{LLZ} utilized it to prove boundary $W^{2,p}$ estimates for linear elliptic equations on $C^{1,\alpha}$ domains.

We illustrate our idea as follows. Let $\{Q_k\}_{k=1}^{\infty}$ be Whitney decomposition of $\Omega_1$ (Suppose $0\in\partial\Omega$ and denote $\Omega_r=\Omega\cap B_r$) and  $\widetilde Q_k=\frac65Q_k$ the $\frac65-$dilation of $Q_k$ with respect to its center. Let $u$ satisfy \eqref{1.1}.
Deduce from  H$\rm{\ddot{o}}$lder's inequality and interior $W^{2,\delta_0}$ estimates that for any $\delta\leq\delta_0$,
\begin{equation*}
\begin{aligned}
||D^2 u||_{L^\delta(Q_k)}^\delta
&\leq Cd_k^{n(1-\frac{\delta}{\delta_0})}||D^2 (u-l)||_{L^{\delta_0}(Q_k)}^\delta\\
&\leq C \left(d_k^{n-2\delta}||u-l||_{L^{\infty}(\widetilde Q_k)}^\delta
+d_k^{n-\frac{\delta n}{p}}||f||_{L^{p}(\widetilde Q_k)}^\delta
\right)
\end{aligned}
\end{equation*}
for some universal constant $C$ and any affine function $l$, where $d_k$ denotes the diameter of $Q_k$.
By $C^{1,\alpha_0}$ estimates up to the boundary, we can take $l$ such that
$$|u(x)-l(x)|\leq  C\dist(x,\partial\Omega)^{1+\alpha_0},\ \forall x\in \Omega_{1/2}$$
for some positive constants $C$ and proper $0<\alpha_0<1$. It then follows from Whitney decomposition that
$$||u-l||^\delta_{L^{\infty}(\widetilde Q_k)}\leq Cd_k^{(1+\alpha_0)\delta}.$$
Combining the above estimates, we obtain
\begin{equation*}
\begin{aligned}
||D^2 u||_{L^\delta(Q_k)}^\delta
\leq C \left(d_k^{-\delta+\alpha_0 \delta+n}+d_k^{n-\frac{\delta n}{p}}||f||_{L^{p}(\Omega_1)}^\delta
\right).
\end{aligned}
\end{equation*}
By taking sum on both sides with respect to $k$, we obtain the desired boundary $W^{2,\delta}$ estimate for any $\delta<\min\{1/(1-\alpha_0),p/n\}$ that guarantees $\sum_k(d_k^{-\delta+\alpha_0 \delta+n}+d_k^{n-\frac{\delta n}{p}})$ is convergent.


The paper is organized as follows. In Section 2, Whitney decomposition and its relevant properties are concluded. In Section 3, we demonstrate some basic
estimates for $u\in S(\lambda,\Lambda,f)$ concerning interior $W^{2,\delta_0}$ estimates and pointwise boundary H$\rm{\ddot{o}}$lder gradient estimates. In Section 4, we prove Theorem \ref{t1.1}.

We end this section by introducing some notations.

\bigskip

\noindent\textbf{Notation.} \\
1. $e_i=(0,...,0,1,...,0)=i^{th}$ standard coordinate vector.\\
2. $x'=(x^1,x^2,...,x^{n-1})$ and $x=(x',x^n).$\\
3. $\mathbb{R}^n_+=\{x\in \mathbb{R}^n:x^n>0\}.$\\
4. $B_r(x_0)=\{x\in \mathbb{R}^n: |x-x_0|<r\}$ and $B^+_r (x_0) = B_r(x_0) \cap  \mathbb{R}^n_+$.\\
5. $B_r'=\{x'\in \mathbb{R}^{n-1}: |x'|<r\}$ and $T_r= \{(x',0): x'\in B_r'\}.$\\
6. $\Omega_r (x_0)= \Omega\cap B_r(x_0)$ and $(\partial\Omega)_r (x_0)= \partial\Omega\cap B_r(x_0)$. We omit $x_0$ when $x_0=0.$\\
7. $\diam E$ = $\sup\{|x-y|,x,y\in E\}$,\  $\forall E\subset\mathbb{R}^n.$\\
8. $\dist(E,F)$ = $\inf\{|x-y|,x\in E\ ,y\in F\}$, \ $\forall E,F\subset\mathbb{R}^n.$\\

\section{Whitney decomposition}
In what follows, by a cube we mean a closed cube in $\mathbb{R}^n$, with sides parallel to the axes. We say two such cubes are disjoint if their interiors are disjoint.

\bigskip

\begin{lemma}\label{l2.1}(Whitney decomposition)
Let $\Omega$ be a non-empty open set in $\mathbb{R}^n$. Then there exist two sequences of cubes $Q_k$(called the Whitney cubes of $\Omega$) and $\widetilde Q_k=\frac65{Q_k}$ ($\frac 65-$dilation of $Q_k$ with respect to center of $Q_k$) such that

(i) $\Omega=\bigcup_{k=1}^{\infty}Q_{k}=\bigcup_{k=1}^{\infty}\widetilde Q_k$;

(ii) The $Q_k$ are mutually disjoint;

(iii) $d_k\leq \mathrm{dist}\ (Q_k,\partial\Omega)\leq 4d_k$, where $d_k=\diam Q_k$;

(iv) Each point of $\Omega$ is contained in at most $12^n$ of the cubes $\widetilde Q_k$.
\end{lemma}

For the proof of the above lemma, we refer to Theorem 1 and Proposition 3 in Section VI.1 in \cite{ST}.
~\\

From now on, we make the following assumption on $\Omega$:
~\\
$\textbf{(A)}$\ $0\in \partial\Omega$ and there exists a continuous function $\varphi: B_1'\rightarrow\mathbb{R}$  such that
\begin{equation*}\label{domain}
\Omega_1=\{x^n>\varphi(x'),|x|<1\}\ \ \ \mathrm{and}\ \ \ (\partial\Omega)_1=\{x^n=\varphi(x'),|x|<1\}.
\end{equation*}
Recall that $\Omega_1=\Omega\cap B_1$ and $(\partial\Omega)_1=\partial\Omega\cap B_1$.
~\\

Let $\{Q_k\}_{k=1}^{\infty}$ be Whitney decomposition of $\Omega_{1}$, $d_k=\diam Q_k$ and ${\widetilde Q_k}=\frac 65 Q_k$.
The following lemmas were first proved in \cite{LLZ} and we give their proofs here for completeness.

\bigskip

\begin{lemma}\label{lf'}
Suppose that $\Omega$ satisfies Assumption $(A)$. Then
\begin{equation}\label{1/4}
\Omega_{1/{12}}\subset\bigcup_{\widetilde Q_k\subset \Omega_{1/{4}}} Q_k.
\end{equation}
\end{lemma}

\begin{proof}
If not, there exist a point $x\in \Omega_{1/12}$ and a cube $Q_k$ such that $x\in Q_k$ but $\widetilde Q_k\not\subset \Omega_{1/4}$. It follows that there exists a point  $y\in \widetilde Q_k$ with $|y|\geq 1/4$. Then we deduce from Lemma \ref{l2.1} (iii) that
$$\dist (Q_k, \partial \Omega_{1})\geq \diam Q_k=5\diam \widetilde Q_k/6\geq{5}(|y|-|x|)/6\geq 5/36>1/12.$$
On the other hand, since $x\in Q_k\cap\Omega_{1/12}$ and $0\in \partial\Omega$,
$$\dist (Q_k, \partial \Omega_{1})\leq|x|\leq 1/12.$$
Thus we get a contradiction.
\end{proof}

\bigskip

\begin{lemma}\label{lf}
Suppose that $\Omega$ satisfies Assumption $(A)$ with $\varphi\in C^{0,1}(B_1')$. If $q>n-1$, then
\begin{equation}\label{sum}
\sum\limits_{\widetilde{Q}_k\subset \Omega_{1/4}} d_k^q\leq C,
\end{equation}
where $C$ depends on $n$, $q$ and $||\varphi||_{C^{0,1}(B_1')}$.
\end{lemma}

\begin{proof}

Set
\begin{equation}\label{fs}
\mathcal{F}^{s}=\bigcup_{k}\ \{Q_k:2^{-s-1}< d_k\leq 2^{-s},\ \widetilde{Q}_k\subset \Omega_{1/4}\},\ s=1,2,....
\end{equation}
For any  $Q_k\in \mathcal{F}^{s}$, since $\mathcal{F}^{s}\subset \Omega_{1/4}$, there exists
$y_k\in (\partial\Omega)_{1/2}$ such that
$$\dist(Q_k,y_k)=\dist(Q_k,(\partial\Omega)_1)=\dist(Q_k,\partial\Omega_1)\leq 4d_k\leq2^{-s+2},$$
where we have used Lemma \ref{l2.1} (iii). (Recall that $(\partial\Omega)_r=\partial\Omega\cap B_r$ for $r>0$.)
It follows that
$$\dist(x,(\partial\Omega)_1)\leq d_k+\dist(Q_k,y_k)
\leq 2^{-s}+2^{-s+2}\leq 2^{-s+3}$$
for any $x\in Q_k$ and $Q_k\in \mathcal{F}^{s}$, which implies
\begin{equation}\label{ss}
\mathcal{F}^{s}\subset\Omega_{1/4}\cap\{\dist(x,(\partial\Omega)_1)\leq 2^{-s+3}\}.
\end{equation}

By Assumption $(A)$, we have
\begin{equation*}
\begin{aligned}
&\Omega _{1/4}\cap\{\dist(x,(\partial\Omega)_1)\leq 2^{-s+3}\}\\
&\mbox{}\hskip0.3cm \subset\{|x'|\leq 1/4,\varphi(x')\leq x^n\leq \varphi(x')+(||\varphi||_{C^{0,1}(B_1')}+1)2^{-s+3}\}.
\end{aligned}
\end{equation*}
Since
\begin{equation*}
|\{|x'|\leq 1/4,\varphi(x')\leq x^n\leq \varphi(x')+(||\varphi||_{C^{0,1}(B_1')}+1)2^{-s+3}\}|\leq C2^{-s},
\end{equation*}
we have
$$|\Omega _{1/4}\cap\{\dist(x,(\partial\Omega)_1)\leq 2^{-s+3}\}| \leq C2^{-s},$$
where $C$ depends on $n$ and $||\varphi||_{C^{0,1}(B_1')}$.
Therefore, by \eqref{ss},
\begin{equation}\label{ss'}
|\mathcal{F}^{s}|\leq C{2^{-s}}.
\end{equation}

Observe that
$$\bigcup_{\widetilde{Q}_k\subset \Omega_{1/4}}Q_k
=\bigcup_{s=1}^{\infty}\bigcup_{Q_k\in \mathcal{F}^{s}} Q_k.$$
If $q>n-1$, we
derive from \eqref{fs} and \eqref{ss'} that
$$\begin{aligned}
\displaystyle\sum\limits_{\widetilde{Q}_k\subset \Omega_{1/4}} d_k^q
&\leq\displaystyle\sum\limits_{s=1}^{\infty}\left\{
\displaystyle\sum\limits_{Q_k\in \mathcal{F}^{s}}
\left(d_k^{q-n}\cdot d_k^{n}\right)\right\}
\leq C \displaystyle\sum\limits_{s=1}^{\infty}\left\{2^{-s(q-n)}
\cdot\displaystyle\sum\limits_{Q_k\in \mathcal{F}^{s}}{d_k^{n}} \right\}\\
\\&\leq C \displaystyle\sum\limits_{s=1}^{\infty}2^{-s(q-n)}|\mathcal{F}^{s}|\leq
C \displaystyle\sum\limits_{s=1}^{\infty}2^{-s(q-n+1)}\leq C,
\end{aligned}$$
where $C$ depends on $n$, $q$ and $||\varphi||_{C^{0,1}(B_1')}$.
\end{proof}

\begin{re}
Lemma \ref{lf} is obvious as $q\geq n$ since
$$\sum\limits_{\widetilde{Q}_k\subset \Omega_{1/4}} d_k^q\leq\sum\limits_{\widetilde{Q}_k\subset \Omega_{1/4}} d_k^n\leq C_n|\Omega_{1/4}|$$
for some dimensional constant $C_n$.
\end{re}

\section{Preliminary estimates}

 The following lemma concerns interior $W^{2,\delta_0}$ regularity for $u\in S(\lambda,\Lambda,f)$ and $f\in L^n$. We refer to \cite{CC}, \cite{CCKS} and \cite{L} for its  proof.

\bigskip

\begin{thm}\label{t1.0}
There exists a positive constant $\delta_0$ depending on $n$, $\lambda$ and $\Lambda$ such that if
$$u\in S(\lambda,\Lambda,f) \ \ \mathrm{in}  \ B_1
\ \ \mathrm{and} \ \ f\in L^n(B_1),$$
then $u\in W^{2,\delta_0}(B_{1/2})$ and
\begin{equation}\label{1.5}
||u ||_{ W^{2,\delta_0}(B_{1/2})}\leq C\left(||u ||_{L^{\infty}(B_1)}+ ||f||_{L^n(B_1)}\right),
\end{equation}
where $C$ depends on $n$, $\lambda$ and $\Lambda$.
\end{thm}

The following lemma concerns pointwise boundary  $C^{1,\bar\alpha}$ estimates for $u\in S(\lambda,\Lambda,0)$ on flat boundaries with zero boundary values. We owe it to Krylov\cite{K} and see  more details in Lian and Zhang \cite{LZ}.

\bigskip

\begin{lemma}\label{l3.1}
Let $u$ satisfy
\begin{equation*}
\left\{
\begin{aligned} &u\in S(\lambda,\Lambda,0)   &\mathrm{in} \ \   &B_1^+;\\
&u=0   &\mathrm{on}  \ \  &T_1.
\end{aligned}
\right.
\end{equation*}
Then there exists a constant $0<\bar \alpha<1$ depending on $n,\lambda$ and $\Lambda$ such that $u$ is $C^{1,\bar \alpha}$ at $0$, i.e., there exists a constant $a$ such that
$$|u(x)-ax^n|\leq C|x|^{1+\bar\alpha}||u||_{L^\infty(B_1^+)},\ \forall x\in B_{1/2}^+$$
and
$$|Du(0)|=|a|\leq C||u||_{L^\infty(B_1^+)},$$
where $C$ depends on $n$, $\lambda$ and $\Lambda$.
\end{lemma}

Based on Lemma \ref{l3.1} and compactness method, Lian and Zhang \cite{LZ} proved the following lemma concerning pointwise boundary  $C^{1,\alpha}$ estimates for $u\in S(\lambda,\Lambda,f)$  on  $C^{1,\alpha}$ domains with $0<\alpha<\bar\alpha$.

\bigskip

\begin{lemma}\label{l3.2}
Let $\bar\alpha$ be as in  Lemma \ref{l3.1} and  $0<\alpha<\bar\alpha$.
Suppose that $\Omega$ satisfies Assumption $(A)$  with $\varphi\in C^{1,\alpha}(B_1')$ and $u$ satisfies
\begin{equation}\label{4.1}
\left\{
\begin{aligned} &u\in S(\lambda,\Lambda,f)   &\mathrm{in} \ \   &\ \Omega_1;\\
&u=g   &\mathrm{on}  \ \  &(\partial\Omega)_1,
\end{aligned}
\right.
\end{equation}
where $g\in C^{1,\alpha}((\partial\Omega)_1)$ and $f\in L^n(\Omega_1)$ satisfying for some constant $K_f$,
\begin{equation}\label{3.3}
||f||_{L^n(\Omega_r)}\leq K_f r^{\alpha}\ \ \mathrm{for\ any}\  r\in(0,1].
\end{equation}
Then $u$ is $C^{1,\alpha}$ at 0, i.e., there exists an affine function $l$ such that
$$|u(x)-l(x)|\leq C|x|^{1+\alpha}
(||u||_{L^\infty(\Omega_1)}+K_f+||g||_{C^{1,\alpha}((\partial\Omega)_1)}),\ \forall x\in \Omega_{r_0}$$
and
$$|Du(0)|=|Dl|\leq C(||u||_{L^\infty(\Omega_1)}+K_f+||g||_{C^{1,\alpha}((\partial\Omega)_1)}),$$
where $C$ and $r_0$ depends on $n$, $\lambda$,  $\Lambda$, $\alpha$, $\bar\alpha$ and $||\varphi||_{C^{1,\alpha}(B_1')}$.
\end{lemma}

If $f\in L^p(\Omega_1)$ for $n<p<\infty$, by H\"{o}lder's inequality, there exists a constant $C_n$ depending only on $n$ such that
\begin{equation*}
||f||_{L^n(\Omega_r)}\leq C_n r^{1-n/p} ||f||_{L^p(\Omega_{r})}\ \ \mathrm{for\ any}\  r\in(0,1],
\end{equation*}
which implies $f$ satisfies \eqref{3.3} with $\alpha=1-n/p$. Then we have the following corollary:
\bigskip

\begin{cor}\label{l3.3}
Under the hypotheses of Lemma \ref{l3.2} with \eqref{3.3} replaced by
$f\in L^p(\Omega_1)$ for $n<p<\infty$.
Set $\alpha_0=\min\{\alpha,1-n/p\}$. Then
$u$ is $C^{1,\alpha_0}$ at 0, i.e.,
there exists an affine function $l_{x_0}$ such that
$$|u(x)-l(x)|\leq C|x|^{1+\alpha_0}
(||u||_{L^\infty(\Omega_1)}+||f||_{L^p(\Omega_1)}+||g||_{C^{1,\alpha}((\partial\Omega)_1)}), \ \forall x\in \Omega_{r_0}$$
and
$$|Du(0)|=|Dl|\leq C(||u||_{L^\infty(\Omega_1)}+||f||_{L^p(\Omega_1)}+||g||_{C^{1,\alpha}((\partial\Omega)_1)}),$$
where $C$ and $r_0$ depends on $n$, $\lambda$, $\Lambda$, $\alpha$, $\bar\alpha$ and $||\varphi||_{C^{1,\alpha}(B_1')}$.
\end{cor}

\section{Proof of  Theorem \ref{t1.1}}

Theorem \ref{t1.1} follows easily from the following Theorem \ref{t3.1} which is proved by interior $W^{2,\delta_0}$  estimates (Theorem \ref{t1.0}), boundary  $C^{1,\alpha_0}$ estimates (Corollary \ref{l3.3}) and  Whitney decomposition (Lemma \ref{l2.1}, \ref{lf'} and \ref{lf}).

\bigskip

\begin{thm}\label{t3.1}
Let $\bar\alpha$ be as in Lemma \ref{l3.1} and $0<\alpha<\bar\alpha$. Suppose that $\Omega$ satisfies Assumption $(A)$  with $\varphi\in C^{1,\alpha}(B_1')$ and $u$ satisfies
\begin{equation*}
\left\{
\begin{aligned} &u\in S(\lambda,\Lambda,f)   &\mathrm{in} \ \   &\ \Omega_1;\\
&u=g   &\mathrm{on}  \ \  &(\partial\Omega)_1,
\end{aligned}
\right.
\end{equation*}
where
$$f\in L^p(\Omega_1)\ \ \mathrm{for} \ \  n<p<\infty \ \ \mathrm{and} \ \  g\in C^{1,\alpha}((\partial\Omega)_1).$$
Let $\delta_0$ be as in Theorem \ref{t1.0} and $\alpha_0=\min\{\alpha,1-n/p\}$. Then for any
\begin{equation}\label{d}
\delta\leq\delta_0 \ \ \mathrm{and}\ \ \delta <1/(1-\alpha_0),
\end{equation}
we have
\begin{equation}\label{u}
||u||_{ W^{2,\delta}(\Omega_{1/12})}\leq C(||u||_{L^{\infty}(\Omega_1)}+||f||_{L^p(\Omega_1)}+ ||g||_{C^{1,\alpha}((\partial\Omega)_1)}),
\end{equation}
where $C$ depends on $n,\lambda,\Lambda,\alpha,\bar\alpha,\delta,p$ and $||\varphi||_{C^{1,\alpha}(B_1')}$.
\end{thm}
\begin{proof}
In the following, we write for simplicity
$$\mathcal{H}=||u||_{L^\infty(\Omega_{1})}+||f||_{L^p(\Omega_1)}
+||g||_{C^{1,\alpha}((\partial\Omega)_{1})}.$$

Let $\{Q_k\}_{k=1}^{\infty}$ be Whitney decomposition of $\Omega_{1}$, $d_k=\diam Q_k$ and ${\widetilde Q_k}=\frac 65 Q_k$.
For any $\widetilde Q_k\subset \Omega_{1/4}$, let $y_k\in (\partial\Omega)_{1/2}$ and $\tilde x_k\in \partial\widetilde Q_k$ such that
$$|\tilde x_k-y_k|=\dist(\widetilde Q_k, \partial\Omega_1)<\dist(Q_k, \partial\Omega_1)\leq 4d_k,$$
where Lemma \ref{l2.1} (iii) is used in the last inequality.
Consequently, we see that
$$ |x-y_k|\leq|x-\tilde x_k|+|\tilde x_k-y_k|\leq 5d_k,\ \ \forall x\in\widetilde Q_k.$$
By  Corollary \ref{l3.3}, $u$ is $C^{1,\alpha_0}$ at $y_k$ and then there exists an affine function $l_{y_k}$(written by $l$ for simplicity in the following) such that
\begin{equation}\label{3.4}
|u(x)-l(x)|\leq C |x-y_k|^{1+\alpha_0}\mathcal{H}\leq  C d_k^{1+\alpha_0}\mathcal{H},\ \forall x\in \widetilde Q_k,
\end{equation}
where $C$ depends on $n,\lambda,\Lambda,\alpha,\bar\alpha$ and $||\varphi||_{C^{1,\alpha}(B_1')}$.

Since
$u-l\in S(\lambda,\Lambda,f)$ in $\widetilde Q_k$, we deduce
from interior $W^{2,\delta_0}$ estimate, \eqref{3.4} and
H$\rm{\ddot{o}}$lder's inequality  that
\begin{equation*}
\begin{aligned}
\displaystyle\int_{Q_k}|D^2 (u-l)|^{\delta_0}dx
&\leq C\left\{
d_k^{n-2{\delta_0}}||u-l||^{\delta_0}_{L^\infty(\widetilde Q_k)}
+d_k^{n-{\delta_0}}
\left(\displaystyle\int_{\widetilde Q_k}|f|^{n}dx\right)^{{\delta_0}/n}
\right\}\\
&\leq C\left\{d_k^{n-(1-\alpha_0){\delta_0}}\mathcal{H}^{\delta_0} +d_k^{n-{\delta_0} n/p}
\left(\displaystyle\int_{\widetilde  Q_k}|f|^{p}dx\right)^{{\delta_0}/p}
\right\},
\end{aligned}
\end{equation*}
where $C$ depends on $n,\lambda,\Lambda,\alpha,\bar\alpha,p$ and $||\varphi||_{C^{1,\alpha}(B_1')}$.

For any $\delta\leq {\delta_0},$  we deduce from H$\rm{\ddot{o}}$lder's inequality and the above estimate that
\begin{equation*}
\begin{aligned}
\displaystyle\int_{Q_k}|D^2 u|^{\delta}dx
&\leq Cd_k^{n(1-\delta/{\delta_0})}
\left(\displaystyle\int_{Q_k}|D^2 (u-l)|^{\delta_0}dx\right)^{\delta/{\delta_0}}\\
&\leq C\left\{d_k^{n-(1-\alpha_0){\delta}}\mathcal{H}^{\delta}+
d_k^{n-{\delta} n/p}
\left(\displaystyle\int_{\widetilde  Q_k}|f|^{p}dx\right)^{{\delta}/p}
\right\}.
\end{aligned}
\end{equation*}
Taking sum on both sides with respect to $k$, we obtain
\begin{equation*}
\begin{aligned}
\displaystyle\sum\limits_{\widetilde Q_k\subset \Omega_{1/4}}
\displaystyle\int_{Q_k}|D^2 u|^{\delta}dx
\leq C\left\{\mathcal{H}^{\delta}\displaystyle\sum\limits_{\widetilde Q_k\subset \Omega_{1/4}}
d_k^{n-(1-\alpha_0){\delta}}
+||f||_{L^p(\Omega_1)}^\delta\displaystyle\sum\limits_{\widetilde Q_k\subset \Omega_{1/4}}d_k^{n-{\delta} n/p}
\right\}.
\end{aligned}
\end{equation*}

Since $\alpha_0=\min\{\alpha,1-n/p\}$ and $\delta<1/(1-\alpha_0)$, we have
$$n-{\delta} n/p\geq n-(1-\alpha_0)\delta>n-1.$$
Consequently, by Lemma \ref{lf'} and Lemma \ref{lf}, for any $\delta\leq\delta_0$ and $\delta<1/(1-\alpha_0),$
\begin{equation*}
\begin{aligned}
\displaystyle\int_{\Omega_{1/12}}|D^2 u|^{\delta}dx
&\leq \displaystyle\sum\limits_{\widetilde Q_k\subset \Omega_{1/4}}
\displaystyle\int_{Q_k}|D^2 u|^{\delta}dx\\
&\leq C\left(||u||_{L^\infty(\Omega_{1})}^\delta+||f||_{L^p(\Omega _1)}^\delta+||g||_{C^{1,\alpha}((\partial\Omega)_{1})}^\delta\right),
\end{aligned}
\end{equation*}
where $C$ depends on $n,\lambda,\Lambda,\alpha,\bar\alpha,\delta,p$ and $||\varphi||_{C^{1,\alpha}(B_1')}$.
The desired estimate \eqref{u} then follows easily from the above estimate.
\end{proof}

Theorem \ref{t1.1} then follows from Theorem \ref{t3.1}.

\bigskip

\noindent\emph{\bf{Proof of Theorem \ref{t1.1}.}}

\bigskip

Since $T\subset\partial\Omega$ is of $C^{1,\alpha}$, for any $x\in T$, there exist $r_x>0$ and $\varphi_x\in C^{1,\alpha}(B'_{r_x}(x'))$ such that
$$\Omega_{r_x}(x)=\{x^n>\varphi_x(x')\}\cap B_{r_x}(x)
 \ \mathrm{and}\ \ (\partial\Omega)_{r_x}(x)=\{x^n=\varphi_x(x')\}\cap B_{r_x}(x).$$
Since $\bigcup_{x\in T}B_{{r_x}/12}(x)$ cover $T\cap \bar\Omega'$ and $T\cap \bar\Omega'$ is compact, we choose a finite subset $B_{r_i/12}(x_i),\ i=1,2,...,N,$ of $\{B_{{r_x}/12}(x):x\in T\}$ that still covers $T\cap \bar\Omega'$.

Let $\delta_0$ be as in Theorem \ref{t1.0}, $\bar\alpha$ be as in Lemma \ref{l3.1} and $\alpha_0=\min\{\bar\alpha,\alpha,1-n/p\}$. Utilizing the scaled version of Theorem \ref{t3.1} in $\Omega_{r_i}(x_i)=\Omega\cap B_{r_i}(x_i)$, we obtain that for any $\delta\leq\delta_0$ and
$\delta<1/(1-\alpha_0)$,
\begin{equation*}
\begin{array}{ccc}
r_i^{-n/\delta-2}||u||_{ L^{\delta}(\Omega_{r_i/12}(x_i))}
+r_i^{-n/\delta-1}||D u||_{ L^{\delta}(\Omega_{r_i/12}(x_i))}
+r_i^{-n/\delta}
||D^2u||_{ L^{\delta}(\Omega_{r_i/12}(x_i))}\\
\leq C_i\left(
r_i^{-2}||u||_{L^{\infty}(\Omega_{r_i}(x_i))}
 +r_i^{-n/p}||f||_{L^{p}(\Omega_{r_i}(x_i))}
+r_i^{-2}||g||_{C^{1,\alpha}((\partial\Omega)_{r_i}(x_i))}\right)
\end{array}
\end{equation*}
for $i=1,2,...,N$, where $C_i$ depends on $n,\lambda,\Lambda,\alpha,\bar\alpha,\delta, p$ and $||\varphi_{x_i}||_{C^{1,\alpha}(B_{r_i}'(x_i'))}$.

Finally, we sum the above inequalities along with the interior $W^{2,\delta_0}$ estimate, to find $u\in W^{2,\delta}(\Omega')$  for $\Omega'\subset\subset\Omega\cup T$, with the desired estimate \eqref{1.2}.
\qquad\qquad\qquad $\Box$

\section*{Acknowledgement}
This work is supported by NSFC 12071365.




\bibliographystyle{elsarticle-num}



\end{document}